\documentclass{article}
\newcommand{\proof}{\noindent {\bf Proof: }}

\newtheorem{theorem}{Theorem}
\newtheorem{statement}{Statement}
\newtheorem{corollary}{Corollary}
\newtheorem{lemma}{Lemma}
\newtheorem{defi}{Definition}

\def\qed{\hfill $\Box$}

\usepackage[]{amssymb}
\usepackage{graphicx}

\begin{document}
\title{Semi-indefinite-inner-product and generalized Minkowski spaces}
\author{\'A.G.Horv\'ath\\ Department of Geometry, \\
Budapest University of Technology and Economics,\\
H-1521 Budapest,\\
Hungary}
\date{Nov. 3, 2008}

\maketitle

\begin{abstract}
In this paper we parallelly build up the theories of normed linear
spaces and of linear spaces with indefinite metric,  called also
Minkowski spaces for finite dimensions in the literature.

In the first part of this paper we collect the common properties of the semi- and indefinite-inner-products and define the
semi-indefinite-inner-product and the corresponding structure, the semi-indefinite-inner-product space. We give  a generalized concept of Minkowski
space embedded in a semi-indefinite-inner-product space using the concept of a new product, that contains the classical cases as special ones.

In the second part of this paper we investigate the real, finite dimensional generalized Minkowski space and its sphere of radius $i$. We prove that
it can be regarded as a so-called Minkowski-Finsler space and if it is homogeneous one with respect to linear isometries, then the Minkowski-Finsler
distance its points can be determined by the Minkowski-product.
\end{abstract}



\section{Introduction}

\subsection{Notation and Terminology}

\begin{description}

\item[concepts without definition:] real and complex vector spaces, basis, dimension, direct sum of subspaces, linear and bilinear mapping,
quadratic forms, inner (scalar) product, hyperboloid, ellipsoid, hyperbolic space and hyperbolic metric, kernel and rank of a linear mapping.
\item[i.p.:] Inner ( or scalar product) of a vector space.
\item[s.i.p.:] Semi-inner-product (See Definition 1)
\item[continuous s.i.p.:] The definition can be seen after Definition 1.
\item[differentiable s.i.p.:] See Definition 3.
\item[i.i.p.:] Indefinite-inner-product (See Definition 4)
\item[s.i.i.p:] Semi-indefinite-inner-product (See Definition 6).
\item[Minkowski product:] See Definition 7.
\item[generalized Minkowski space:] See Definition 7.
\item[generalized space-time model:] Finite dimensional, real, generalized Minkowski space with one dimensional time-like orthogonal direct components.
\item[positive (resp. negative) subspace:] It is a subspace in an i.i.p. space in which all vectors have positive (resp. negative) scalar square.
\item[neutral or isotropic subspace:] See Definition 5.
\item[Auerbach basis:] The corresponding definition with respect to a finite dimensional real normed space can be seen before Theorem 8.
\item[hypersurface:] The definition in a generalized Minkowski space can be seen before Lemma 3.
\item[tangent vector, tangent hyperplane:] These definitions can be seen before Lemma 3.
\item[Minkowski-Finsler space:] See Definition 15.
\item[$\mathbb{C}$, $\mathbb{R}$, $\mathbb{R}^{n}$, $S^n$:] The complex line, the real line, the $n$-dimensional
    Euclidean space and the $n$-dimensional unit sphere,
    respectively.
\item[$\dim(V)$]: The dimension of the vector space $V$.
\item[$x\bot y$:] The notion of (non-symmetric) property of orthogonality. We consider it in the meaning "y is orthogonal to x".
\item[{$[\cdot,\cdot]$}:] The notion of scalar product and all its generalization.
\item[{$[\cdot,\cdot]^-$}:] The notion of s.i.p. corresponding to a generalized Minkowski space.
\item[{$[\cdot,\cdot]^+$}:] The notion of Minkowski product of a generalized Minkowski space.
\item[{$[x,\cdot]'_z(y)$}:] The derivative map of an s.i.p. in its second argument into the direction of $z$ at the point $(x,y)$. See Definition 3.
\item[$\|\cdot\|'_x(y)$,$\|\cdot\|''_{x,z}(y)$:]The derivative of the norm in the direction of $x$ at the point $y$ and the second derivative of the norm in the directions $x$ and $z$ at the point $y$.
\item[$<\{\cdot\}>$:] The linear hull of a set.
\item[$\Re{\{\cdot\}}$,$\Im{\{\cdot\}}$:] The real and imagine part of a complex number, respectively.
\item[$T_v$:] The tangent space of a Minkowskian hypersurface in its point $v$.
\item[$\mathcal{S},\mathcal{T},\mathcal{L}$:]The set of space-like, time-like and light-like vectors respectively.
\item[$S$,$T$:]Space-like and time-like orthogonal direct components of a generalized Minkowski space, respectively.
\item[$\{e_1,\ldots,e_k,e_{k+1},\ldots,e_n\}$:] An Auerbach basis of a generalized Minkowski space with $\{e_1,\ldots,e_k\}\subset S$ and $\{e_{k+1},\ldots,e_n\}\subset T$, respectively.
\item[$H$,$H^+$:] Are the sphere of radius $i$ and its upper sheet, respectively.
\end{description}

\subsection{Completion of the preliminaries }

In this introduction we recall some important moments from the long history of the theme of this paper. We complete these some observations are needed for our investigations.

\subsubsection{Semi-inner-product spaces}

A generalization of the inner product and the inner product spaces (briefly i.p spaces) raised by G.Lumer in \cite{lumer}.
\begin{defi}[\cite{lumer}]
The semi-inner-product (s.i.p) on a complex vector space $V$ is a complex function $[x,y]:V\times V\longrightarrow \mathbb{C}$ with the following properties:
\begin{description}

\item[s1]: $[x+y,z]=[x,z]+[y,z]$
\item[s2]: $[\lambda x,y]=\lambda[x,y]$ \mbox{ for every } $\lambda \in \mathbb{C}$
\item[s3]: $[x,x]>0$ \mbox{ when } $x\not =0$
\item[s4]: $|[x,y]|^2\leq [x,x][y,y]$
\end{description}
A vector space $V$ with a s.i.p. is a s.i.p. space.
\end{defi}

He proved that a s.i.p space is a normed vector space with norm $\|x\|=\sqrt{[x,x]}$ on the other hand every normed vector space can be represented
as a s.i.p. space. In \cite{giles} J.R.Giles showed that a homogeneity property:

\begin{description}

\item[s5]: $[x,\lambda y]=\bar{\lambda}[x,y]$ for all complex $\lambda $,

\end{description}
can be imposed, and all normed vector spaces can be represented as s.i.p. spaces with this property. Giles also introduced the concept of {\bf continuous s.i.p. space} as s.i.p. space having the additional property:

\begin{description}
\item[s6]: For every unit vectors $x,y \in S$, $\Re\{[y,x+\lambda y]\}\rightarrow\Re\{[y,x]\}$ for all real $\lambda\rightarrow 0$.
\end{description}
The space is uniformly continuous if the above limit is approached uniformly for all points $x,y$ of the unit sphere $S$.

A characterization of the continuous s.i.p. space is based on the differentiable property of the space.
\begin{defi}[\cite{giles}]
A normed space is G\^{a}teaux differentiable if for all $x,y$ elements of its unit sphere and real value $\lambda$,
$$
\lim\limits_{\lambda \rightarrow 0}\frac{\|x+\lambda y\|-\|x\|}{\lambda} \mbox{ exists.}
$$
A normed vector space is uniformly Fr\`{e}chet differentiable if this limit is approached uniformly for the pair $x,y$ points of the unit sphere.
\end{defi}

Giles in \cite{giles} proved that
\begin{theorem}[\cite{giles}]
An s.i.p. space is continuous (uniformly continuous) s.i.p. space if and only if the norm is G\^{a}teaux (uniformly Fr\`{e}chet) differentiable.
\end{theorem}

In the second part of this paper we need a stronger condition on differentiability of the s.i.p. space, therefore we define the differentiable s.i.p. as follows:

\begin{defi}
A differentiable s.i.p. space is an s.i.p. space where the s.i.p. has the additional property:

\noindent {\bf s6':} For every three vectors x,y,z and real $\lambda $
$$
[x,\cdot]'_z(y):=\lim\limits_{\lambda \rightarrow 0}\frac{\Re\{[x,y+\lambda z]\}-\Re\{[x,y]\}}{\lambda } \mbox{ does exist.}
$$
We say that the s.i.p. space is continuously differentiable, if the above limit as a function of $y$, is continuous.
\end{defi}

First we note that the equality $\Im\{[x,y]\}=\Re\{[-ix,y]\}$ with the above property guaranties the existence and continuity of the complex limit:
$$
\lim\limits_{\lambda \rightarrow 0}\frac{[x,y+\lambda z]-[x,y]}{\lambda }.
$$
Analogously to the theorem of Giles (see Theorem 3 in \cite{giles}) we connect this definition to the differentiability properties of the norm function generated by the s.i.p..

\begin{theorem}
An s.i.p. space is a (continuously) differentiable s.i.p. space if and only if the norm is two times (continuously) G\^{a}teaux differentiable.
\end{theorem}

Before the proof of the theorem we introduce a notion on G\^{a}teaux derivative of the norm. Let
$$
\|\cdot\|'_x(y):=\lim\limits_{\lambda \rightarrow 0}\frac{\|y+\lambda x\|-\|y\|}{\lambda },
$$
be the derivative of the norm in the direction of $x$ at the point $y$. Similarly we use the notation:
$$
\|\cdot\|''_{x,z}(y):=\lim\limits_{\lambda \rightarrow 0}\frac{\|\cdot\|'_x(y+\lambda z)-\|\cdot\|'_x(y)}{\lambda }
$$
which is the second derivative of the norm in the directions $x$ and $z$ at the point $y$.
We need the following useful lemma going back, with different notation to McShane \cite{mcshane} or Lumer \cite{lumer2}:

\begin{lemma}[\cite{lumer2}]
If E is any s.i.p. space, $x,y\in E$, then:
$$
\|y\| (\|\cdot\|'_x(y))^-\leq \Re\{[x,y]\}\leq \|y\| (\|\cdot\|'_x(y))^+
$$
where $(\|\cdot\|'_x(y))^-$ and $(\|\cdot\|'_x(y))^+$ denotes the left hand and right hand derivatives with respect to the real variable $\lambda $. In particular, if the norm is differentiable, then
$$
[x,y]= \|y\| \{(\|\cdot\|'_x(y))+\|\cdot\|'_{-ix}(y)\}.
$$
\end{lemma}

\proof[of Theorem 2]
To determine the derivative of the s.i.p. assume that the norm is differentiable twice. Then by the Lemma 1, above:
$$
\frac{\Re\{[x,y+\lambda z]\}-\Re\{[x,y]\}}{\lambda }=\frac{\|y+\lambda z\|(\|\cdot\|'_x(y+\lambda z))-\|y\| (\|\cdot\|'_x(y))}{\lambda }=
$$
$$
=\frac{\|y\|\|y+\lambda z\|(\|\cdot\|'_x(y+\lambda z))-\|y\|^2 (\|\cdot\|'_x(y))}{\lambda \|y\|}\geq
$$
$$
\geq \frac{|[y+\lambda z,y]|(\|\cdot\|'_x(y+\lambda z))-\|y\|^2 (\|\cdot\|'_x(y))}{\lambda \|y\|},
$$
 where we have assumed that the sign of $\frac{\|\cdot\|'_x(y+\lambda z)}{\lambda}$ is positive. Since the derivative of the norm is continuous this follows from the assumption that $\frac{\|\cdot\|'_x(y)}{\lambda}$ is positive. Considering the latter condition we get:
$$
\frac{\Re\{[x,y+\lambda z]\}-\Re\{[x,y]\}}{\lambda }\geq
$$
$$
\geq \|y\|^2 \frac{\|\cdot\|'_x(y+\lambda z)-(\|\cdot\|'_x(y))}{\lambda \|y\|}+\frac{\Re[z,y]}{\|y\|}\|\cdot\|'_x(y+\lambda z).
$$
On the other hand,
$$
\frac{\|y+\lambda z\|(\|\cdot\|'_x(y+\lambda z))-\|y\| (\|\cdot\|'_x(y))}{\lambda }\leq
$$
$$
\leq \frac{\|y+\lambda z\|^2(\|\cdot\|'_x(y+\lambda z))-|[y,y+\lambda z]| (\|\cdot\|'_x(y))}{\lambda \|y+\lambda z\|}=
$$
$$
=\frac{\|y+\lambda z\|^2(\|\cdot\|'_x(y+\lambda z))- (\|\cdot\|'_x(y))}{\lambda \|y+\lambda z\|}+\lambda \Re{[z,y+\lambda z]}\frac{(\|\cdot\|'_x(y))}{\lambda \|y+\lambda z\|}.
$$
Analogously, if $\frac{\|\cdot\|'_x(y)}{\lambda}$ is negative, then both of the above inequalities are revers, we get that the limit
$$
\lim\limits_{\lambda \mapsto 0}\frac{\Re\{[x,y+\lambda z]\}-\Re\{[x,y]\}}{\lambda } \mbox{ exists,}
$$
and equals to
$$
\|y\|(\|\cdot\|''_{x,z}(y))+\frac{\Re[x,y]\Re[z,y]}{\|y\|^2}.
$$
Here we note that in the case  when $\frac{\|\cdot\|'_x(y)}{\lambda}=0$ also there does exist a neighborhood in which the sign of the function $\frac{\|\cdot\|'_x(y+\lambda z)}{\lambda}$ is constant. Thus we don't have to investigate this case by itself.
Conversely, consider the fraction:
$$
\|y\|\frac{\|\cdot\|'_x(y+\lambda z)-(\|\cdot\|'_x(y))}{\lambda}.
$$
We assume now that the s.i.p. is differentiable implying that it is continuous, too. The norm is differentiable by the theorem of Giles. Using again Lemma 1 and  assuming that $\frac{\Re[x,y]}{\lambda}> 0$ we have:
$$
\|y\|\frac{\|\cdot\|'_x(y+\lambda z)-(\|\cdot\|'_x(y))}{\lambda}=\frac{\Re[x,y+\lambda z]\|y\|-\Re[x,y]\|y+\lambda z\|}{\lambda \|y+\lambda z\|}=
$$
$$
=\frac{\Re[x,y+\lambda z]\|y\|^2-\Re[x,y]\|y+\lambda z\|\|y\|}{\lambda \|y\|\|y+\lambda z\|}\leq \frac{\Re[x,y+\lambda z]\|y\|^2-\Re[x,y]|[y+\lambda z,y]|}{\lambda \|y\|\|y+\lambda z\|}=
$$
$$
=\frac{\Re\{[x,y+\lambda z]\}-\Re\{[x,y]\}}{\lambda }\frac{\|y\|}{\|y+\lambda z\|}-\frac{\Re[x,y]\Re[z,y]}{\|y\|\|y+\lambda z\|}.
$$
On the other hand using the continuity of the s.i.p. and our assumption $\frac{\Re[x,y]}{\lambda}> 0$ similarly as above, we also get an inequality:
$$
\|y\|\frac{\|\cdot\|'_x(y+\lambda z)-(\|\cdot\|'_x(y))}{\lambda}\geq
$$
$$
\frac{\Re\{[x,y+\lambda z]\}-\Re\{[x,y]\}}{\lambda} -\frac{\Re[x,y+\lambda z]\Re[z,y+\lambda z]}{\|y+\lambda z\|^2}.
$$
If we reverse the assumption of signs then the direction of inequalities will change, too. A limit argument shows again that the first differential function is differentiable and the connection between the two derivatives is:
$$
\|y\|(\|\cdot\|''_{x,z}(y))= [x,\cdot]'_z(y)-\frac{\Re[x,y]\Re[z,y]}{\|y\|^2}.
$$
\qed

\subsubsection{Further remarks on the theory of s.i.p.}

B.Nath gave a straightforward generalization of a s.i.p. by replacing the Schwartz's inequality by the H\"older's inequality in \cite{nath}. He showed that this kind of generalized s.i.p. space induces a norm by setting $\|x\|=[x,x]^{\frac{1}{p}}$ $1\leq p\leq \infty $, and for every normed space it can be constructed a generalized s.i.p. space. (For $p=2$, this theorem reduces to Theorem 2 of Lumer.) The connection between the Lumer-Giles s.i.p. and the generalized s.i.p. of Nath is simple. The s.i.p. $[x,y]$ for every $p's$ defines a generalized s.i.p. by the equality:
$$
\widehat{[x,y]}=[y,y]^{\frac{p-2}{p}}[x,y].
$$
The s.i.p. holds the homogeneity property of Giles if and only if the Nath's generalized s.i.p. satisfies the $p-1$-homogeneity property:

\begin{description}

\item[s5'']: $\widehat{[x,\lambda y]}=\bar{\lambda}|\lambda|^{p-2}\widehat{[x,y]}$ for all complex $\lambda $.

\end{description}
Thus in this paper we will concentrate only the original version of the s.i.p..

From geometric point of view if $K$ is a $0$-symmetric, bounded,
convex body in the Euclidean $n$-space $\mathbb{R}^n$ (with a fixed
origin O) then it defines a norm whose unit ball is $K$ itself (see
\cite{l-g}). Such a space is called Minkowski normed space. The main
results in this topic are collected in the surveys
\cite{martini-swanepoel 1}, \cite{martini-swanepoel 2} and
\cite{martini}. In fact, the norm is a continuous function which is
considered (in geometric terminology as in \cite{l-g}) as a gauge
function. Combining  this with the result of Lumer and Giles we get
that a Minkowski normed space can be represented as an s.i.p space.
The metric (the so-called Minkowski metric), i.e. the distance of
two points induced by this norm, is invariant with respect to the
translations of the space.

\subsubsection{Indefinit inner product spaces}

Another concept of  Minkowski space raised by H.Minkowski and used
by in theoretical physic and differential geometry based upon the
concept of the indefinite inner product. (See e.g. \cite{gohberg}.)

\begin{defi}[\cite{gohberg}]
The indefinite inner product (i.i.p.) on a complex vector space $V$
is a complex function $[x,y]:V\times V\longrightarrow \mathbb{C}$
with the following properties:
\begin{description}

\item[i1]: $[x+y,z]=[x,z]+[y,z]$
\item[i2]: $[\lambda x,y]=\lambda[x,y]$ \mbox{ for every } $\lambda \in \mathbb{C}$
\item[i3]: $[x,y]=\overline{[y,x]}$ \mbox{ for every } $x,y\in V$
\item[i4]: $[x,y]=0$ \mbox{ for every } $y\in V$ then $x=0$.
\end{description}
A vector space $V$ with an i.i.p. is an i.i.p. space.
\end{defi}

We recall, that a subspace in an i.i.p. space is positive (nonnegative) if all of its nonzero vectors have
positive (nonnegative) scalar squares. The classification of subspaces with respect to the positivity property in an i.i.p. space is also an
interesting question. First we pass now to the class of subspaces which are peculiar to i.i.p. spaces and have no analogous in the spaces with a
definite inner product.
\begin{defi}[\cite{gohberg}]
A subspace $N$ in $V$ is called neutral if $[v,v]=0$ for all $v\in N$.
\end{defi}
In view of the identity
$$
[x,y]=\frac{1}{4}\{[x+y,x+y]+i[x+iy,x+iy]-[x-y,x-y]-i[x-iy,x-iy]\}
$$
a subspace $N$ is neutral in an i.i.p. space if and only if
$[u,v]=0$ for all $u,v\in N$. Observe also that a neutral subspace
is nonpositive and nonnegative in the same time, and is necessarily
degenerate. Therefore the following statement can be proved:
\begin{theorem}[\cite{gohberg}]
An nonnegative (resp. nonpositive) subspace is a direct sum of a
positive (resp. negative) subspace and a neutral subspace.
\end{theorem}
We note that the decomposition of a nonnegative subspace $U$ into a
direct sum to a positive and a neutral component is not unique, in
general. However, the dimension of the positive summand is uniquely
determined.

The standard  mathematical model of the space-time is a four
dimensional i.i.p. space with signature $(+,+,+,-)$. This is also
called Minkowski space in the literature.

\subsection{Results}

In the first  part of this paper we introduce the concept of semi-indefinite-inner-product (s.i.i.p.) and the generalized notation of Minkowski
space. We also define the concept of orthogonality of such spaces. (Section 2.)

In the second part  we give the definition of the Minkowski-Finsler space in a generalized space-time model. This construction goes similarly to the
definition of a Riemannian manifold (e.g. geometric Minkowski space or hyperbolic space) by embedding into an i.i.p. space. (Section 3.)

We prove only those statements whose proof can not to be found in the literature. (These are: Statement1, Theorems 2, 7-11, 13-15 and Lemmas 2-4,
respectively.) The author uses the known statements without proof and gives references to them.

\section{Unification and geometrization}
\subsection{Semi-indefinite-inner product spaces}

In this section  let {\bf s1}, {\bf s2}, {\bf s3}, {\bf s4}, be the
four defining properties of a s.i.p and {\bf s5} be the homogeneity
property of the second argument imposed by Giles, respectively.
(Namely, {\bf s1} is the additivity property of the first argument,
{\bf s2} is the homogeneity property of the first argument, {\bf s3}
meaning the positivity of the function, {\bf s4} is the
Cauchy-Schwartz inequality.)

On the other hand  {\bf i1}={\bf s1}, {\bf i2}={\bf s2}, {\bf i3} is the antisymmetry property and {\bf i4} is the nondegeneracy property of the
product, respectively. It is easy to see that {\bf s1}, {\bf s2}, {\bf s3}, {\bf s5} imply {\bf i4} and if $N$ is a positive (negative) subspace of a
i.i.p. space then {\bf s4} holds on $N$. In the following definition we connect the concepts of s.i.p. and i.i.p..

\begin{defi}
The semi-indefinite-inner-product (s.i.i.p.) on a complex vector space $V$ is a complex function $[x,y]:V\times V\longrightarrow \mathbb{C}$ with the following properties:
\begin{description}

\item[1] $[x+y,z]=[x,z]+[y,z]$ (additivity in the first argument)
\item[2] $[\lambda x,y]=\lambda[x,y]$ \mbox{ for every } $\lambda \in \mathbb{C}$ (homogeneity in the first argument)
\item[3] $[x,\lambda y]=\overline{\lambda}[x,y]$ \mbox{ for every } $\lambda \in \mathbb{C}$ (homogeneity in the second argument)
\item[4] $[x,x]\in \mathbb{R}$ \mbox{ for every } $x\in V$ (the corresponding quadratic form is real valued)
\item[5] If either $[x,y]=0$ \mbox{ for every } $y\in V$ or $[y,x]=0$ for all $y\in V$ then $x=0$ (nondegeneracy)
\item[6] $|[x,y]|^2\leq [x,x][y,y]$ holds on nonpositive and nonnegative subspaces of V, respectively. (Cauchy-Schwartz inequality is valid on positive and negative subspaces, resp.)
\end{description}
A vector space $V$ with a s.i.i.p. is a s.i.i.p. space.
\end{defi}
Interest in s.i.i.p. spaces depends largely on the example spaces given by the s.i.i.p. space structure.

\noindent {\bf Example 1:} We conclude that an s.i.i.p. space is a
homogeneous s.i.p. space if and only if the property {\bf s3} holds,
too. An s.i.i.p. space is an i.i.p. space if and only if the
s.i.i.p. is an antisymmetric product. In this latter case
$[x,x]=\overline{[x,x]}$ implies {\bf 4}, and the function is
Hermitian linear in its second argument, too. In fact, we have:
$[x,\lambda y+\mu z]=\overline{[\lambda y+\mu
z,x]}=\overline{\lambda} \overline{[y,x]}+\overline{\mu} \overline{[
z,x]}=\overline{\lambda}[x,y]+\overline{\mu}[x,z]$. It is clear that
the classical "Minkowski spaces" can be represented by either a
s.i.p or an i.i.p, so they automatically can be represented as an
s.i.i.p. space, too.

\noindent {\bf Example 2:} Let now $V=<\{e_1,\ldots,e_n\}>$ be a finite dimensional vector space and $C$ be the surface of a cross-polytope defined
by:
$$
C=\cup \{ \mbox{ conv }\{\varepsilon_ie_i | i=1,\ldots,n\} \mbox{ for all choices of } \varepsilon_i=\pm 1\}.
$$
It is clear that for a real vector $v\in C$ there exists at least
one linear functional, and we choose exactly one $v^\star$ of the
dual space holding the property $v^\star(v)=(-1)^k$ where $k$ is the
combinatorial dimension of that combinatorial face $F_v$ of C which
contains the point $v$ in its relative interior. (It is easy to see
that $k+1$ is the cardinality of the nonzero coefficients of the
representation of v.) For $\lambda v\in V$ where $v\in C$ and any
real $\lambda $ (by Giles method) we choose $(\lambda
v)^\star=\lambda v^\star$. Given such a mapping from $V$ into
$V^\star$, it is readily verified that the product
$$
 [u,v]=v^\star(u)
$$
satisfies the properties {\bf 1}-{\bf 4}. {\bf 5} also  holds since
there is no vector $v$ for which $v^\star(v)=0$. Finally, every
two-dimensional subspace has vectors $v$ and $w$ by $v^\star(v)>0$
and $w^\star(w)<0$ there are  neither positive nor negative
subspaces with dimension at least two, implying that property {\bf
6} holds, too.

\noindent {\bf Example 3:} In an arbitrary complex normed linear space $V$ we can define an s.i.i.p. which is a generalization of a representing
s.i.p. of the norm function. Let now $C$ be the unit sphere of the space $V$. By the Hahn-Banach theorem
there exists at least one continuous linear functional, and we choose exactly one such that $\|\widetilde{v}^\star\|=1$ and
$\widetilde{v}^\star(v)=1$. Consider a sign function $\varepsilon (v)$ with value
$\pm 1$ on $C$. If now $\varepsilon ({v})=1$ let denote by $v^\star=\widetilde{v}^\star$  and if $\varepsilon ({v})=-1$ define
$v^\star=-\widetilde{v}^\star$. Finally, homogeneously extract it to $V$
by the equality $(\lambda v)^\star=\overline{\lambda}v^\star$ as in the previously example. Of course for an arbitrary vector $v$ of $V$ the corresponding linear functional satisfies the equalities
$v^\star(v):=\varepsilon (v)\|v\|^2$ and $\|v\|=\|v^\star\|$. Now the function
$$
[u,v]=v^\star(u)
$$
satisfies {\bf 1}-{\bf 5}. If $U$ is a nonnegative subspace then it is positive and we have for all nonzero $u,v\in U$:
$$
|[u,v]|=|v^\star(u)|=\frac{|v^\star(u)|}{\|u\|}\|u\|\leq \|v^\star\|\|u\|=\|v\|\|u\|,
$$
proving {\bf 6}.

\subsection{The generalized Minkowski space}

Before the definition we prove an important lemma.

\begin{lemma}
Let  $(S,[\cdot,\cdot]_S)$ and $(T,-[\cdot,\cdot]_T)$ be two s.i.p.
spaces. Then the function
$[\cdot,\cdot]^-:(S+T)\times(S+T)\longrightarrow \mathbb{C}$ defined
by
$$
[s_1+t_1,s_2+t_2]^-:=[s_1,s_2]-[t_1,t_2]
$$
is an s.i.p. on the vector space $S+T$.
\end{lemma}

\proof The function $[\cdot,\cdot]^-$ is nonnegative, as we can see
from its definition easily. First we prove the linearity in the
first argument. We have:
$$
[\lambda'(s'+t')+\lambda''(s''+t''),s+t]^-=[\lambda's'+\lambda''s'',s]_S-[\lambda't'+\lambda''t'',t]_T=
$$
$$
=\lambda'[s',s]_S+\lambda''[s'',s]_S-\lambda'[t',t]_T-\lambda''[t'',t]_T=\lambda'[s'+t',s+t]^-+\lambda''[s''+t'',s+t]^-.
$$
The homogeneity in the second argument is trivial. In fact,
$$
[s'+t',\lambda(s+t)]^-=[s',\lambda s]_S-[t',\lambda
t]_T=\overline{\lambda}[s'+t',s+t]^-
$$
Finally we check the inequality of Cauchy-Schwartz. Since we have:
$$
|[s_1+t_1,s_2+t_2]^-|^2=[s_1+t_1,s_2+t_2]^-\overline{[s_1+t_1,s_2+t_2]^-}=
$$
$$
=([s_1,s_2]_S-[t_1,t_2]_T)(\overline{[s_1,s_2]_S}-\overline{[t_1,t_2]_T})=
$$
$$
=[s_1,s_2]_S\overline{[s_1,s_2]_S}+[t_1,t_2]_T\overline{[t_1,t_2]_T}+[s_1,s_2]_S(-\overline{[t_1,t_2]_T})
+(-[t_1,t_2]_T)\overline{[s_1,s_2]_S}\leq
$$
$$
\leq
[s_1,s_1]_S[s_2,s_2]_S+[t_1,t_1]_T[t_2,t_2]_T+2\Re\{[s_1,s_2]_S(-\overline{[t_1,t_2]_T})\}\leq
$$
$$
\leq [s_1,s_1]_S[s_2,s_2]_S+[t_1,t_1]_T[t_2,t_2]_T+2|[s_1,s_2]_S|
|[t_1,t_2]_T|\leq
$$
$$
\leq[s_1,s_1]_S[s_2,s_2]_S+[t_1,t_1]_T[t_2,t_2]_T+2\sqrt{[s_1,s_1]_S[s_2,s_2]_S[t_1,t_1]_T[t_2,t_2]_T},
$$
and by the inequality between the arithmetic and geometric means we get that:
$$
[s_1,s_1]_S[s_2,s_2]_S+[t_1,t_1]_T[t_2,t_2]_T+2\sqrt{[s_1,s_1]_S[s_2,s_2]_S[t_1,t_1]_T[t_2,t_2]_T}\leq,
$$
$$
\leq
[s_1,s_1]_S[s_2,s_2]_S+[t_1,t_1]_T[t_2,t_2]_T+[s_1,s_1]_S(-[t_2,t_2]_T+(-[t_1,t_1]_T)[s_2,s_2]_S=
$$
$$
=([s_1,s_1]_S-[t_1,t_1]_T)([s_2,s_2]_S-[t_2,t_2]_T)=[s_1+t_1,s_1+t_1]^-[s_2+t_2,s_2+t_2]^-.
$$
\qed

It is possible that the s.i.i.p. space $V$ is a direct sum of its
two subspaces where one of them is positive and the other one is a
negative. Then we have two other structures on $V$, ( by Lemma 2) an s.i.p. structure and a natural third one which we
will call minkowskian  structure. More precisely:

\begin{defi}
Let $(V,[\cdot,\cdot])$ be an s.i.i.p. space. Let $S,T\leq V$ be positive and negative subspaces, where $T$ is a direct complement of $S$ with
respect to $V$. Define a product on $V$ by the equality $[u,v]^+=[s_1+t_1,s_2+t_2]^+=[s_1,s_2]+[t_1,t_2]$, where $s_i\in S$ and $t_i\in T$,
respectively.  Then we say that the pair $(V,[\cdot,\cdot]^+)$ is a generalized Minkowski space with Minkowski product $[\cdot,\cdot]^+$. We also say
that $V$ is a real generalized Minkowski space if it is a real vector space and the s.i.i.p. is a real valued function.
\end{defi}

\begin{remark}

\begin{enumerate}

\item The Minkowski product defined by the above equality satisfies the properties {\bf 1}-{\bf 5} of the s.i.i.p.. But in general property {\bf 6}
does not hold. To see this define a s.i.i.p. space on the following manner: 

Consider a 2-dimensional $L^\infty$ space $S$ of the embedding three
dimensional Euclidean space $E^3$. Choose an orthonormed basis $\{e_1,e_2,e_3\}$ of $E^3$ for which $e_1,e_2\in S$ and give a s.i.p. associated to the
$L^\infty$ norm as follows:
$$
[x_1e_1+x_2e_2,y_1e_1+y_2e_2]_S:=
$$
$$
=x_1y_1\lim\limits_{p\rightarrow\infty}\frac{1}{\left(1+\left(\frac{y_2}{y_1}\right)^p\right)^\frac{p-2}{p}}+
x_2y_2\lim\limits_{p\rightarrow\infty}\frac{1}{\left(1+\left(\frac{y_1}{y_2}\right)^p\right)^\frac{p-2}{p}}.
$$
By Lemma 2 the function
$$
[x_1e_1+x_2e_2+x_3e_3,y_1e_1+y_2e_2+y_3e_3]^-:=[x_1e_1+x_2e_2,y_1e_1+y_2e_2]_S+x_3y_3
$$
is an s.i.p. on $E^3$ associated to the norm
$$
\sqrt{[x_1e_1+x_2e_2+x_3e_3,x_1e_1+x_2e_2+x_3e_3]^-}:=\sqrt{\max\{|x_1|,|x_2|\}^2+x_3^2}.
$$
By the method of Example 3 consider such a sign function for which $\varepsilon(v)$ is equal to $1$ if $v$ is in $S\cap C$ and is equal to $-1$ if
$v=e_3$ holds. ($C$ denotes the unit sphere as in the previous examples.) This sign function determine an s.i.i.p. $[\cdot,\cdot]$ and thus a Minkowski product $[\cdot,\cdot]^+$, for which the generated
square root function is:
$$
f(v):=\sqrt{[x_1e_1+x_2e_2+x_3e_3,x_1e_1+x_2e_2+x_3e_3]^+}=\sqrt{\max\{|x_1|,|x_2|\}^2-x_3^2}.
$$
\begin{figure}
\centerline{\includegraphics[scale=0.5]{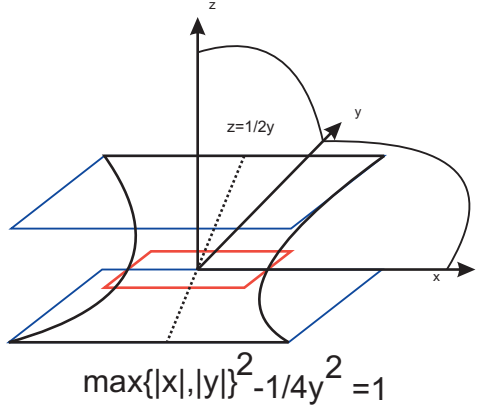}}
\caption{\label{fig1}The unit sphere of a positive subspace of the example in Remark 1}
\end{figure}
As it can be seen easily the plane $x_3=\alpha x_2$ for $0<\alpha<1$ is positive subspace with respect to the Minkowski product but its unit ball is not convex. (See Fig.1) But $f(v)$ homogeneous thus it is not subadditive. Since Cauchy-Scwartz inequality implies subadditivity, this inequality also
false in this positive subspace.

\item The real generalized Minkowski space is a geometrical Minkowski space if it is finite dimensional and the s.i.i.p. is an s.i.p..
(Also implying that its negative component is trivial.) Its Minkowski functional is generated by the norm mapping $\|v\|:v\longmapsto \sqrt{[v,v]}$.
The unit ball of this space is $\{v | [v,v]=1\}$.

\item  The finite dimensional real generalized Minkowski space is a pseudo-Euclidean space if the s.i.i.p is an i.i.p, a space-time model if it is
pseudo-Euclidean and its negative direct component has dimension 1. Its signature of corresponds to the dimensions of $S$ and $T$.

\item By Lemma 2 the s.i.p. $\sqrt{[v,v]^-}$ is a norm function on $V$ which can give an embedding space for a generalized Minkowski
space. This situation is analogous with the situation when a pseudo-Euclidean space is obtained from an Euclidean space by the action of an i.i.p..(See Fig.2.)
\begin{figure}
\centerline{\includegraphics[scale=0.5]{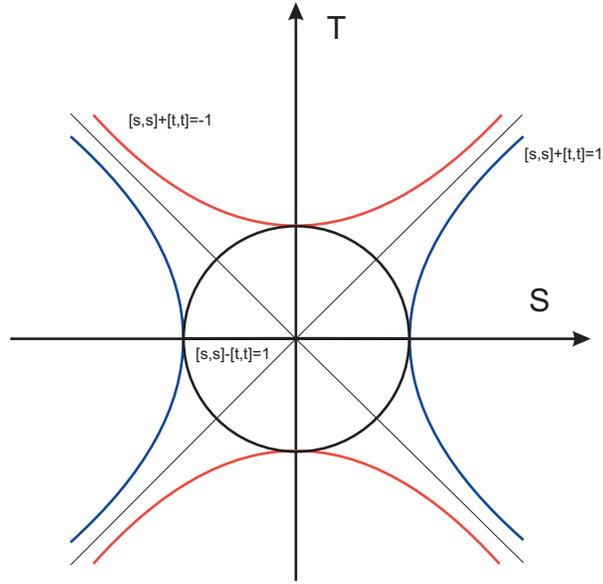}}
\caption{\label{fig2}The real and imaginary unit spheres in dimension two.}
\end{figure}
\end{enumerate}
\end{remark}

\subsection{Further examples for non-trivial s.i.i.p and Minkowski spaces}
\subsubsection{$C^2$ normsquare function  and the associated s.i.i.p. space}

In this section (by Theorem 4) we give a method to construct s.i.i.p spaces with more differentiable property.

A $C^2$ Minkowski space is an $n$-dimensional affine space with metric $d(x,y)=F(y-x)$ where $F$ is the (Minkowskian)
norm function of the associated vector space, where
\begin{description}
\item[n1] $F(x)>0$ for $x\neq 0$
\item[n2] $F(\lambda x)=|\lambda|F(x)$, for all real $\lambda $,
\item[n3] $F(x+y)\leq F(x)+F(y)$, equality holds for $x,y\neq 0$ if and only if $y=\lambda x$ for some real $\lambda >0$
\item[n4] $F(x)$ is of class $C^2$ in each of its $n$ arguments, the components of vector $x$.
\end{description}
Giles in his paper \cite{giles} proved that there is a natural form for an s.i.p. in the associated vector space for which it is a uniform s.i.p. space. The importance of uniform s.i.p. spaces is based on the fact that in such a space the representation theory of Riesz holds and its dual space is also uniform. Now we define a similar class of s.i.i.p. spaces associated to the concept of $C^2$ normsquare function.
\begin{defi}
Consider $\mathbb{R}^n$ as a real vector space $V$ and let $G:V\longrightarrow \mathbb{R}$ be a function on it. If it satisfies the following two properties:
\begin{description}
\item[pn1] $G(\lambda x)=\lambda ^2G(x)$ for real $\lambda $
\item[pn2] If $G|_W\geq (\leq ) 0$ on a subspace $W$ of $V$ then for the positive function $\sqrt{G|_W}$ ($\sqrt{-G|_W}$) holds the convexity property {\bf [n3]},
\end{description}
then we say that $G$ is a normsquare function on $V$. If we also require for $G$ the differentiability property {\bf [n4]}, then we say that the normsquare function is a $C^2$ one.
\end{defi}
It is easy to see that the square of a norm function is a normsquare function, and every i.i.p. defines a normsquare function by $G(x)=[x,x]$. For $C^2$ normsquares we have:
\begin{theorem}
If G is a $C^2$ normsquare function on the real vector space $V$ then there is an associated s.i.i.p. which gives uniform s.i.p. structures on positive (resp. negative) subspaces of $V$.
\end{theorem}
\proof
From the derivatives of a homogeneous function of order 2, for $G$ we have,
$$
DG|_{\lambda x}x=2\lambda G(x) \mbox{ and } x^TD^2G|_{\lambda x}x=2G(x),
$$
where $D(G)_{x}$ means the totally (Frechet) derivative of the function $G$ at the point $x$.
Substituting into these formula $\lambda =1$ we get:
$$
G(x)=\frac{1}{2}x^TD^2(G|_{x})x=\frac{1}{2}DG|_{x}x.
$$
Let the associated s.i.i.p. be defined by the equality:
$$
[x,y]=\frac{1}{2}x^TD^2G|_{y}y.
$$
It is easy to see that this function satisfies properties {\bf 1},{\bf 2},{\bf 4},{\bf 5} of a s.i.i.p.. Property {\bf 3} follows from the fact that $D^2G|_{\lambda x}$ is independent from the value of $\lambda $.
Finally property {\bf 6} is established from the imposed differentiability property and the convexity property {\bf pn2} as follows:
It is clear that the function $\sqrt{G|_W}:W\longrightarrow \mathbb{R}^+$ is a homogeneous $C^2$ function. So
we have:
$$
D\sqrt{G}|_{x}x=\sqrt{G}(x) \mbox{ and } x^TD^2\sqrt{G}|_{x}x=0.
$$
From the identity
$$
D^2G|_{x}=2(\sqrt{G}(x)D^2(\sqrt{G}|_{x})+D\sqrt{G}|_{x}^TD\sqrt{G}|_{x}),
$$
we get that
$$
\frac{1}{2}D^2G|_{y}y=\sqrt{G}(y)D^2\sqrt{G}|_{y}y+D\sqrt{G}|_{y}^TD\sqrt{G}|_{y}y=
$$
$$
=D\sqrt{G}|_{y}^T\sqrt{G}(y)=\sqrt{G}(y)D\sqrt{G}|_{y}^T.
$$
Thus
$$
|[x,y]|=|\frac{1}{2}x^TD^2G|_{y}y|=|x^T\sqrt{G}(y)D\sqrt{G}|_{y}^T|=
$$
$$
=\sqrt{G}(y)|x^TD\sqrt{G}|_{y}^T|=\sqrt{G}(y)|D\sqrt{G}|_{y}x|.
$$
But we have by the second Mean Value Theorem that
$$
\sqrt{G}(x)=\sqrt{G}(y)+D\sqrt{G}|_{y}(x-y)+(x-y)^TD^2\sqrt{G}|_{y+\theta(x-y)}(x-y),
$$
where $0<\theta<1$. Since for a convex $C^2$ function the last summand is non-negative we have that
$$
D\sqrt{G}|_{y}x\leq \sqrt{G}(x),
$$
implying that
$$
|D\sqrt{G}|_{y}x|\leq \sqrt{G}(x).
$$
Thus
$$
|[x,y]|\leq \sqrt{G}(y)\sqrt{G}(x)=\sqrt{[x,x][y,y]},
$$
as we stated. Now the last statement is a consequence of Giles results in \cite{giles}.
\qed

If we have a normed vector space with an associated symmetric, bilinear function then the positive semi-definiteness of the function implies the
Cauchy-Schwartz inequality. If the associated function linear in its first argument and  homogeneous in its second one, the semi-definiteness
property alone does not implies the Cauchy-Schwartz inequality as we can see in the following example.

\noindent{\bf Example 4:} Let $V$ be a two-dimensional vector space with the Euclidean norm:
$$
\|(x,y)^T\|:=\sqrt{x^2+y^2},
$$
where the coordinates can be computed with respect to a fixed orthonormed basis. It is easy to see that an associated product is:
$$
[u_1,u_2]= (x_1\cdot x_2+ 2y_1\cdot y_2)\frac{x^2_2+y^2_2}{x^2_2+2y^2_2},
$$
where $u_i=(x_i,y_i)^T$. This function linear in its first argument homogeneous in its second one, and associated to the norm. On the other hand for
$u_1=(1,2)^T$ and $u_2=(1,1)^T$,
$$
[(1,2)^T,(1,1)^T]=\frac{10}{3}>\sqrt{10}=\sqrt{[(1,2)^T,(1,2)^T]}\sqrt{[(1,1)^T,(1,1)^T]}
$$
gives a counterexample for the Cauchy-Schwartz inequality. The reason of this situation that the norm of the linear functional associated to the
first argument of the product and the fixed vector $u_2$ is greater then the norm of the vector $u_2$.

\subsubsection{Minkowski spaces generated by $L_p$ norms}

\begin{figure}[htbp]
\centerline{\includegraphics[scale=0.5]{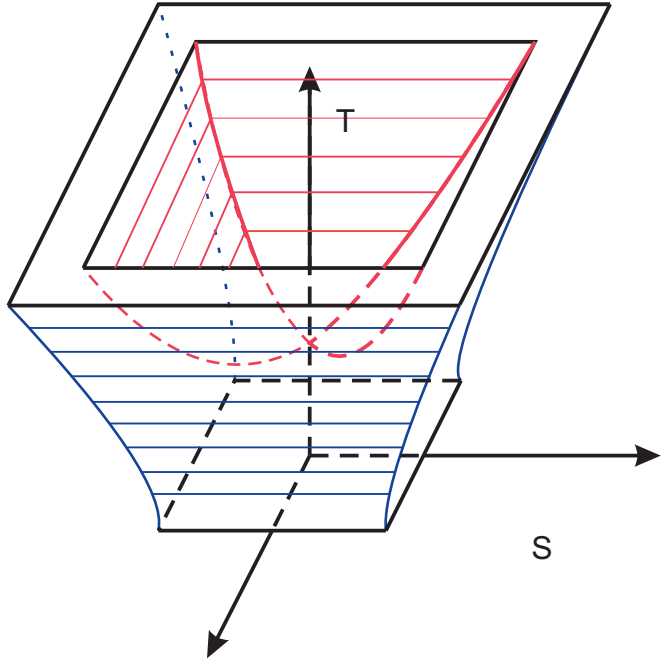}}
\caption{\label{fig3}The case of the norm $L_\infty$.}
\end{figure}

Giles in his paper \cite{giles} gave an associated s.i.p. for $L_p$ spaces. Using the method of our Example 3 we can define s.i.i.p spaces based on
$L_p$ structure. Let $(S,[\cdot,\cdot]_S)$ be the s.i.p. spaces where $S$ is the real Banach space $L_{p_1}(X, \mathcal{S}, \mu)$ and $T$ is the
real Banach space $L_{p_2}(Y, \mathcal{S'}, \nu)$, respectively. If $1<p_1,p_2\leq\infty$ then these spaces can be readily expressed as a uniform s.i.p.
space with s.i.p. defined by
$$
[s_1,s_2]_S=\frac{1}{\|s_2\|^{p_1-2}_{p_1}}\int_X s_1|s_2|^{p_1-1}\mbox{sgn }(s_2) d\mu,
$$
and
$$
[t_1,t_2]_T=\frac{1}{\|t_2\|^{p_2-2}_{p_2}}\int_Y t_1|t_2|^{p_2-1}\mbox{sgn } (t_2) d\nu,
$$
respectively. Consider the real vector space $S+T$ with the s.i.p.
$$
[u,v]^-:=[s_1,s_2]_S+[t_1,t_2]_T.
$$
This is also a uniform s.i.p. space since in Lemma 2 we proved that it is a s.i.p. space and
$$
|[z,x]-[z,y]|=|([s_3,s_1]_S-[s_3,s_2]_S)+([t_3,t_1]_T-[t_3,t_2]_T)|\leq
$$
$$
\leq|[s_3,s_1]_S-[s_2,s_1]_S|+|[t_3,t_1]_T-[t_2,t_1]_T|\leq
$$
$$
\leq 2(p_1-1)\|s_1-s_2\|_{p_1}+2(p_2-2)\|t_1-t_2\|_{p_2},
$$
implying that the space is uniformly continuous. It has been established that such spaces are uniformly convex (\cite{clarkson}, p. 403). By the
method of Example 3 we can define an s.i.i.p spaces on $S+T$ such that the subspace $S$ be positive and $T$ be negative one, and a Minkowski space by
the Minkowski product:
$$
[u,v]^+:=[s_1,s_2]_S-[t_1,t_2]_T,
$$
respectively. (On Fig.3 shows the case when $\dim S=\dim T+1=2$ and the norm of $S$ is $L_\infty$.)

It is easy to see that by this method from every two normed spaces $S$ and $T$ can be done generalized Minkowski space, of course the smoothness
property of it basically determined by the analogous properties of $S$ and $T$.

\subsection{Orthogonality}
\subsubsection{Orthogonality in a normed linear space}

We now investigate an interesting classical problem, the problems of orthogonality in a normed linear space. There are several definitions of orthogonality in a normed linear space which is not an inner product space (i.p. space), but you can not find a concept which is more natural than the others. First we note
that the generalization of the usual i.p. concept of orthogonality
is not unique, that is every concept of orthogonality in s.i.p.
space can be regarded reasonable if it gives back the usual
orthogonality in i.p. sense. Thus we have a lot of possibilities to
define orthogonality. Some of these can be found in the papers
\cite{alonso1}, \cite{alonso2},\cite{james}, \cite{partington},
\cite{diminnie}, \cite{milicic}, \cite{gahler}, \cite{shoja}. Now we
recall the most important concepts.

Let $(V,\|\cdot\|)$ be a normed space and $x,y\in V$. Denote by $x\bot y$ the expression "$y$ orthogonal to "x".
\begin{description}
\item[R] $x\bot y$ iff $\|x+\lambda y\|=\|x-\lambda y\|$ for any $\lambda \in R$ (Roberts, 1934);
\item[B] $x\bot y$ iff $\|x\|\leq\|x+\lambda y\|$ for any $\lambda \in R$ (Birkhoff, 1935);
\item[B-J] $x\bot y$ iff $\|x\|\leq\|x+\lambda y\|$ for any $\lambda \in C$ (Birkhoff-James, 1935);
\item[D] $x\bot y$ iff $\sup\{f(x)g(y)-f(y)g(x) : f,g\in S_{V^\star}\}=\|x\|\|y\|$ where $S_{V^\star}$ denotes the unit sphere of the dual space $V^\star$ (Diminnie, 1983);
\item[A] $x\bot y$ iff $\|x\|\|y\|=0$ or they are linearly independent and such that the four sectors defined by $x$ and $y$ in the unit ball of the plane generated by them (identified to $R^2$) are of the same area (Area, 1984);
\item[I]  $x\bot y$ iff $\|x+y\|=\|x-y\|$ (James' isosceles, 1945);
\item[P]  $x\bot y$ iff $\|x\|^2+\|y\|^2=\|x-y\|^2$ (Pythagorean, 1945);
\item[S]  $x\bot y$ iff $\|x\|\|y\|=0$ or $\|x\|^{-1}x$ and  $\|y\|^{-1}y$ are isosceles orthogonal to each other (Singer or unitary isosceles, 1957);
\item[C]  $x\bot y$ iff $\sum\limits_{i=1}^m\alpha_i\|\beta_i x+\gamma_i y\|^2=0$, where $\alpha_i, \beta_i, \gamma_i$ are real numbers such that
    $$
    \sum\limits_{i=1}^m\alpha_i\beta_i ^2=\sum\limits_{i=1}^m\alpha_i\gamma_i ^2=0, \sum\limits_{i=1}^m\alpha_i\beta_i\gamma_i=1
    $$
    (Carlsson, 1962);
\end{description}

We notice that the {\bf I} and {\bf P} orthogonalities are special cases of the {\bf C} orthogonality. All of the {\bf C} orthogonalities there are an unitary variation as we saw for the isosceles one. The unitary orthogonalities are homogeneous, respectively. This means that if $x\bot y$ then for every pairs of $\lambda,\nu\in \mathbb{C}(\mathbb{R})$ $\lambda x\bot\nu y$ also holds.  It is clear that every orthogonality relation satisfies nondegeneracy ($\lambda x\bot \nu x$ iff either $\lambda x=0$ or $\nu x=0$), simplification (if $x\bot y$, then $\lambda x\bot \lambda y$ for all $\lambda\in \mathbb{R}$) and continuity (if $(x_n)\bot(y_n)\subset V$ such that $x_n\bot y_n$ for every $n\in N$, $x_n\rightarrow x$ and $y_n\rightarrow y$, then $x\bot y$) properties, respectively. On the other hand there are several notions of orthogonality which do not satisfy the symmetric properties. In this case, it is important to distinguish the concept of existence (if $x,y\in V$, then there exists $a\in R$ such that $x\bot (ax+y)$ ) and the concept of additivity (if $x\bot y$ and $x\bot z$ then $x\bot (y+z)$ ) to the left and to the right.

If now we consider the theory of s.i.p in the sense of Lumer-Giles, we have a natural concept of orthogonality. For the unified terminology we change the original calling of Giles and we say that,

\begin{defi}
The vector $y$ is orthogonal to the vector $x$ if $[y,x]=0$.
\end{defi}
Since s.i.p. is neither antisymmetric in the complex case nor symmetric in the real one, this definition of orthogonality is not symmetric in general.

Giles proved that in a continuous s.i.p. space $x$ is orthogonal to $y$ in the sense of the s.i.p. if and only if $x$ is orthogonal to $y$ in the sense of {\bf B-J}. We note that the s.i.p. orthogonality implies the {\bf B-J} orthogonality in every normed spaces. Lumer pointed out that a normed linear space can be transformed into a s.i.p. space in an unique way if and only if its unit sphere is smooth, (i.e. there is an unique supporting hyperplane at each point of the unit sphere). In this case the corresponding (unique) s.i.p. holds the homogeneity property {\bf [s5]}. Imposing the additivity property of the second argument
\begin{description}
\item[s5']: For every $x,y,z\in V$ $[x,y+z]=[x,y]+[x,z]$
\end{description}
the s.i.p. will be a bilinear function. But if the s.i.p. is the unique representation of a given norm and it is bilinear, then it is antisymmetric (resp. symmetric) in the complex (resp. real) case. In fact, define the function $[x,y]':V\times V\longrightarrow \mathbb{C}$ by the equality: $[x,y]'=\overline{[y,x]}$. The properties {\bf s1}, {\bf s2}, {\bf s3}, {\bf s5} trivially hold for this function and the inequality
$$
[x,y]'\overline{[x,y]'}=\overline{[y,x]}[y,x]\leq [y,y][x,x]=[y,y]'[x,x]'
$$
shows the veracity of {\bf s4}. By the unicity of the s.i.p. $[\cdot,\cdot]'$ is equal to the original one, so the s.i.p. is antisymmetric (resp. symmetric), consequently the space is a Hilbert space. (It is an i.p. space.) Summarizing we can say that an unique s.i.p. which is not an i.p. is not additive in its second argument.

\begin{remark}
An orthogonality relation which arises from a s.i.p. representation
of the norm necessarily has the homogeneity property, therefore the
types of orthogonalities in sense of {\bf C} could not be
represented by a s.i.p.. Exactly, it can be proved, that {\bf C}
orthogonality is homogeneous if and only if the space is an i.p.
one. There are a lot of known results and open problems, connecting
with the investigation of the orthogonalities, but as we saw, the
s.i.p. orthogonality of pair of vectors essentially coincides with
their {\bf B-J} orthogonality in the represented normed space. In
this paper we would like to generalize s.i.p. so we have to
concentrate only to the {\bf B-J} orthogonality.

\end{remark}

Another interesting problem is the orthogonality of subspaces. It is clear, that each of the orthogonality relations gives an orthogonality for the subspaces of $V$.
\begin{defi}
Let $X,Y\leq V$ be two subspaces. We say that $X$ is orthogonal to $Y$ if for every pairs of vectors $x\in X$ and $y\in Y$   $x$ is orthogonal to $y$.
\end{defi}
It can be proved that the strongest subspace orthogonality criterium is the Pythagorean one.
\begin{statement}
With respect to subspaces the Pythagorean orthogonality implies any other orthogonality relations.
\end{statement}

\proof We will prove that if two one dimensional subspaces are orthogonal to each other in the sense of Pythagorean orthogonality then the subspace spanned by them is an i.p. space. From this it follows, that all of the 2-dimensional sections of the unit ball of the norm spanned by a vector of $X$ and an other vector of $Y$, are ellipses. This implies that every orthogonality relation restricted to such a plane gives the same orthogonal pairs of vectors as the corresponding i.p., and in this way the original pairs of the two lines are orthogonal to each other in this sense, too.

Consider now the two-plane spanned by the examined linearly independent subspaces $V$ and $V'$. In this plane we consider the usual Euclidean norm and the corresponding i.p. A pair of orthogonal (with respect the Pythagorean orthogonality of the original norm) unit vectors of the subspaces $V$ and $V'$ will be denoted by $v$ and $v'$, respectively. Let $u$ and $u'$ are orthogonal unit vectors with respect to the Euclidean norm. The linear mapping $L$ sends $v$ and $v'$ into the vectors $u$ and $u'$, respectively. If $x$ is an arbitrary unit vector with respect to the original norm we have:
$$
1=\|x\|^2=\|\lambda v +\nu v'\|^2=\|\lambda v\|^2+\|\nu v'\|^2=\lambda^2+\nu ^2
$$
and
$$
L(x)=\lambda u +\nu u'.
$$
This implies that the unit circle of the original norm is mapped onto the unit circle of the Euclidean one, by a linear mapping. Thus the unit circle of the examined plane is a conic. On the other hand the unit circle by our assumption is bounded, consequently is an ellipse. That is the norm is originated from an i.p, as we stated.
\qed

\subsubsection{Orthogonality in an i.i.p. space}

In an i.i.p. space there is a natural definition of the orthogonality.
\begin{defi}[\cite{gohberg}]
Let $(V,[\cdot,\cdot])$ be an i.i.p. space and $U$ be any subset of $V$. Define the orthogonal companion of $U$ in $V$ by
$$
U^\bot =\{v\in V | [v,u]=0 \mbox{ for all } u\in U\}.
$$
\end{defi}
Clearly, $U^\bot$ is a subspace in $V$, and we particularly interested in the case when $U$ is itself a subspace of $V$. In the latter case, it is not generally true that $U^\bot$ is a direct complement for $U$. In contrast, it is true that, for any subspace $U$, the sum of the dimensions of the subspaces $U$ and  $U^\bot$ is equal to the dimension of $V$. The exact answer for this problem uses the concept of nondegeneracy of a subspace, it means that the i.i.p. restricted to this subspace is also nondegenerate. The statement is the following one:
\begin{theorem}[\cite{gohberg}]
$U^\bot$ is a direct complement to $U$ in $V$ if and only if $U$ is nondegenerate.
\end{theorem}

In particular, the orthogonal companion of a nondegenerate subspace is again nondegenerate.

In an i.p. space a fundamental role is played by the construction of a mutually orthogonal set of vectors $u_1,\ldots,u_n$ for which each subset $u_1,\ldots,u_k$ ($k\leq n$) spans the same subspace as a subsets of a given linearly independent set. The well-known Gram-Schmidt process is of this kind. Motivated by applications, attention will be confined to sets of vectors $u_1,\ldots,u_n$ for which $[u_i,u_i]\neq 0$ for each $i$. (Such a vector is called nonneutral.) Note first of all that any set of nonneutral vectors which is orthogonal is necessarily linearly independent. This leads to the concept of regular orthogonalization.

A system of vectors $u_1,\ldots,u_n$ which are mutually orthogonal is said to be a {\bf regular orthogonalization} of $v_1,\ldots,v_n$ if it contains only nonneutral vectors with the property:
$$
<\{u_1,\ldots,u_k\}>= <\{v_1,\ldots,v_k\}>, \mbox{ for } k=1,\ldots ,n.
$$
For any system of vectors $\{v_1,\ldots,v_k\}$, the Gram matrix is defined to the $k\times k$ matrix of the pairwise scalar product of the vectors of the system. The basic statement on regular orthogonalization is the following:
\begin{theorem}[\cite{gohberg}]
The system of vectors $\{v_1,\ldots,v_n\}$ admits a regular orthogonalization if and only if the determinant of its Gram matrix is nonzero. This orthogonalization is essentially unique, if we have two such orthogonal system of vectors then their elements distinct only a scalar factor. (With respect to the complex field $\mathbb{C}$.)
\end{theorem}

\subsubsection{Orthogonality in s.i.i.p spaces}

In this section the pair $(V,[\cdot,\cdot])$ represents a s.i.i.p. space, where $V$ is a complex (real) vector space. We define the orthogonality of such a space with the definition analogous to the definition of the orthogonality of an i.i.p. or s.i.p. space:
\begin{defi}
The vector $v$ is orthogonal to the vector $u$ if $[v,u]=0$. If $U$ is a subspace of $V$, define the orthogonal companion of $U$ in $V$ by
$$
U^\bot =\{v\in V | [v,u]=0 \mbox{ for all } u\in U\}.
$$
\end{defi}

We note that as in the i.i.p. case the orthogonal companion is always a subspace of $V$. The following Theorem analogous to Theorem 4 for i.i.p. spaces.

\begin{theorem} Let $V$ be an $n$-dimensional s.i.i.p. space. Then the orthogonal companion of a nonneutral vector $u$ is a subspace having a direct complement of the linear hull of $u$ in $V$. The orthogonal companion of a neutral vector $v$ is a degenerate subspace of dimension $n-1$ containing $v$.
\end{theorem}

\proof First we observe that if the vector $u$ is nonneutral and its subspace $U=<\{u\}>$, then
$$
U^\bot=\{v | [v,\lambda u]=0 \mbox{ for all } \lambda \in \mathbb{C}\}=\{v | [v,u]=0\}.
$$
Thus $U^\bot\cap U=\emptyset$. On the other hand let the transformation $A:V\longrightarrow V$ defined by $A:x\mapsto [x,u]u$. Obviously it is linear, because of the linearity in the first argument of a s.i.i.p.. Its kernel is
$$
Ker A=\{x | [x,u]u=0\}=\{x | [x,u]=0\}=U^\bot ,
$$
and its image is
$$
Im A=\{[x,u]u | \mbox{  } x\in V\}.
$$
Clearly $Im A$ is a subset of $U$. Since it is a subspace and is not
a trivial one (e.g. $[u,u]u\neq 0$ by our assumption) it is equal to
$U$. By the rank theorem on linear mapping we have that the
dimension of $U^\bot$ is $(n-1)$ and $V$ is a direct sum of $U^\bot$
and $U$.

For a neutral vector $v$ the above argument says that the kernel of
$A$ contains $v$, too. Thus we get $<\{v\}>\subset <\{v\}>^\bot$. On
the other hand taking into consideration the nondegeneracy of $V$
$\dim Im A\neq 0$. Thus again $\dim Im A=1$ and $\dim
<\{v\}>^\bot=(n-1)$ as we stated. \qed

\begin{remark} Observe that this proof does not use the property {\bf 6} of the
s.i.i.p.. So this statement true for any concepts of product
satisfying properties {\bf 1}-{\bf 5}. As we saw, the Minkowski product is also such a product. 
\end{remark}

The following theorem will be a common generalization of the theorem
on diameters conjugated to each other in a real, finite dimensional
normed linear space, and Theorem 6 on the existence of an orthogonal
system in an i.i.p. space. A set of $n$ diameters of the unit ball
of an $n$-dimensional real normed space is considered to be a set of
conjugate diameters if their normalized vectors have the following property:
Choosing one of them, each vector in the linear span of
the remaining direction vectors orthogonal to it. An {\bf Auerbach basis} of a normed
space is a set of direction vectors having this property. Any real
normed linear space has at least two Auerbach bases. One is induced
by a cross-polytope inscribed in the unit ball of maximal volume
(\cite{taylor}), and the other by the midpoints of the facets of a
circumscribed parallelotope of minimum volume (\cite{day}). These
two ways of finding Auerbach bases are dual in the sense that if an
Auerbach basis is induced by an inscribed cross-polytope of maximal
volume, then any dual basis is induced by a circumscribed
parallelotope of minimum volume, and vice versa (\cite{knowles}). If
any minimum volume basis and maximum volume basis coincide, then by
a result of Lenz (\cite{lenz}) we have that the space is a real i.p.
space of finite dimension.

In a generalized Minkowski spaces we have an analogous theorem:

\begin{theorem}
In a finite dimensional, real generalized Minkowski space there is a
basis with the Auerbach property. With other words, its vectors are
orthogonal to the $(n-1)$-dimensional subspace spanned by the
remaining ones. For this basis there is a natural number $k$ less or
equal to $n$, for which $\{e_1,\ldots, e_k\}\subset S$ and
$\{e_{k+1},\ldots, e_n\}\subset T$. Finally, this basis also has the
Auerbach property in the s.i.p. space $(V,[\cdot,\cdot]^-)$.
\end{theorem}

\proof Consider an Auerbach basis in $\{e_1,\ldots, e_k\}\subset S$
in the real normed space generated by the s.i.i.p. in $S$ and
another one $\{e_{k+1},\ldots, e_n\}\subset T$ in the other normed
space generated by the negative of the s.i.i.p. on $T$. The union of
these bases is an Auerbach basis for the Minkowski product and the
s.i.p. $[\cdot,\cdot]^-$, respectively. In fact, e.g. the vectors of
the linear hull of ${e_2,\ldots, e_n}$ are orthogonal to $e_1$,
since
$$
[\alpha_2e_2+\cdots+\alpha_k e_k+\beta_{k+1}
e_{k+1}+\cdots+\beta_ne_n,e_1]^+=
$$
$$
=[\alpha_2e_2+\cdots+\alpha_k e_k,e_1]+[\alpha_{k+1}e_{k+1}+\cdots+\alpha_n e_n,0]=0
$$
is valid by the Auerbach property of ${e_1,\ldots, e_k}$. On the
other hand we have the equalities:
$$
[e_i,e_j]^-=[e_i,e_j]=0 \mbox{ for } 1\leq i,j \leq k,
$$
$$
[e_i,e_j]^-=-[e_i,e_j]=0 \mbox{ for } k+1\leq i,j \leq n
$$
and
$$
[e_i,e_j]^-=0 \mbox{ otherwise }.
$$
This proves the last statement of the theorem.
\qed

\begin{corollary}
In a generalized Minkowski space the positive and negative components $S$ and $T$ are orthogonal to each other in the sense of
Pythagorean orthogonality. In fact, for every pair of vectors $s\in S$ and $t\in T$, by definition we have
$[s-t,s-t]^+=[s,s]+[-t,-t]=[s,s]^++[t,t]^+$.
\end{corollary}

\section{Generalized space-time model and its imaginary unit sphere}

In this section we consider a special subset, the imaginary unit sphere of a finite dimensional, real, generalized Minkowski space. (Some steps of
our investigation is valid in a complex generalized Minkowski space, too. If we don't use the attribute "real" then we think about a complex
Minkowski space.) We give a metric on it and thus we will get a structure similar to the hyperboloid model of the hyperbolic space embedding in a
space-time model. A similar building up of the hyperboloid model of the hyperbolic geometry can be found e.g. in \cite{cannon}.

\begin{defi}
Let $V$ be a generalized Minkowski space. Then we call a vector space-like, light-like and time-like if its scalar square is positive, zero or
negative, respectively. Let denote by $\mathcal{S}, \mathcal{L}$ and $\mathcal{T}$ the sets of the space-like, light-like and time-like vectors,
respectively.
\end{defi}

In a finite dimensional, real generalized Minkowski space for which $\dim T=1$ we can characterize geometrically these sets of vectors. Such a space
is called {\bf generalized space-time model}. In this case $\mathcal{T}$ is a union of its two parts,
$$
 \mathcal{T}=\mathcal{T}^+\cup \mathcal{T}^-
$$
where
$$
\mathcal{T}^+=\{t\in \mathcal{T} | \mbox{ where } t=\lambda e_n \mbox{ for } \lambda \geq 0\} \mbox{ and }
$$
$$\mathcal{T}^-=\{t\in \mathcal{T} |
\mbox{ where } t=\lambda e_n \mbox{ for } \lambda \leq 0\}.
$$
\begin{theorem}
Let $V$ be a generalized space-time model. Then $\mathcal{T}$ is an open double cone with boundary $\mathcal{L}$ and the positive part $\mathcal{T}^+$ (resp. negative part $\mathcal{T}^-$) of $\mathcal{T}$ is convex.
\end{theorem}

\proof The conic property immediately follows from the equality:
$$
[\lambda v,\lambda v]^+=\lambda \overline{\lambda }[v,v]^+=|\lambda |^2[v,v]^+.
$$
Consider now the affine subspace of dimension $(n-1)$ which is of the form $U=S+t$, where $t\in T$ arbitrary, but non zero. Then for an element of
$\mathcal{T}\bigcap U$ we have
$$
0\geq [s+t,s+t]^+=[s,s]+[t,t]
$$
and therefore that $[s,s]\leq -[t,t]$. This implies that the above intersection is a convex body on the $(n-1)$-dimensional real vector space $S$. The s.i.i.p. in $S$ induces a norm whose unit ball is a centrally symmetric convex body. So
$\mathcal{T}$ is a double cone and its positive (resp. negative ) part is convex as we stated. For the vectors of its boundary the equality property holds thus these are light-like vectors. Since those vectors of the space for which the inequality does not hold, are space-time vectors, we also get the remaining statement of the theorem.
\qed

\subsection{The imaginary unit sphere $H$.}

We note that if $\dim T> 1$ or the space is complex then the set of time-like vectors can not be divided into two convex components so we have to
consider that our space is a generalized space-time model.
\begin{defi}
The set
$$
H:=\{ v\in V | [v,v]^+=-1\},
$$
is called the imaginary unit sphere.
\end{defi}

As we saw with respect to the embedding real normed linear space
$(V,[\cdot,\cdot]^-)$ (see Lemma 2) $H$ is a generalized two sheets
hyperboloid corresponding the two piece of $\mathcal{T}$,
respectively. Usually we deal only with one sheet of the hyperboloid
or identify the two sheets projectively. In this case the space-time
component $s\in S$ of $v$ determines uniquely the time-like  one
$t\in T$. Let $v\in H$ be arbitrary. Let denote by $T_v$ the set
$v+v^\bot$ where $ v^\bot$ is the orthogonal complement subspace of
$v$ with respect to the s.i.i.p..

\begin{theorem}
The set $T_v$ corresponding to the point $v=s+t\in H$ is a positive (n-1)-dimensional affine subspace of the generalized Minkowski space
$(V,[\cdot,\cdot]^+)$.
\end{theorem}

\proof By the definition of $H$ the component $t$ of $v$ is non-zero. As we saw in the Remark after Theorem 7 if $[v,v]\neq 0$ then $v^\bot$ is an $(n-1)$-dimensional subspace of $V$. Let now $w\in T_v-v$ be an arbitrary vector. We have to prove that if $[v,v]=-1$ and $w$ orthogonal to $v$ then $[w,w]>0$. Let now $w=s'+t'$ and assume that $[t',t']=0$. Then by the definition of $T$ $t'=0$ and thus $[w,w]=[s,s]>0$ holds. Thus we may assume that $[t',t']\neq 0$ and so $t'=\lambda t$. On the other hand we have:
$$
0=[w,v]^+=[s',s]+[t',t].
$$
We can use the Cauchy-Schwartz inequality for the space-time components, so we have:
$$
[s,s][s',s']\geq |[s',s]|^2=|-[t',t]|^2=|\lambda |^2|-[t,t]|^2=|\lambda |^2[t,t]^2.
$$
Since
$$
[s,s][t',t']=\lambda\overline{\lambda}[s,s][t,t]=|\lambda |^2[s,s][t,t],
$$
we get the inequality:
$$
[s,s][w,w]^+=[s,s]([s',s']+[t',t'])\geq |\lambda |^2([t,t]^2+[s,s][t,t]).
$$
By the definition of $H$ we also have,
$$
-1=[v,v]^+=[s,s]+[t,t],
$$
so
$$
[s,s][w,w]^+\geq |\lambda |^2([t,t]^2+(-1-[t,t])[t,t])=-|\lambda |^2[t,t]>0.
$$
Consequently, if $s$ is nonzero then $[w,w]>0$ as we stated.

If now $[s,s]=0$ then $[t,t]=-1$ and $0=[s'+t',t]=[s',t]+[t',t]=[t',t]$ implies that $t'=0$ and $w\in S$. Thus we proved the statement.
\qed

Each of the affine spaces $T_v$ of $H$ can be considered as a semi-metric space, where the semi-metric arises from the Minkowski product restricted
to this positive subspace of $V$. We recall that Minkowski product does not hold the Cauchy-Schwartz inequality thus the corresponding distance function
does not hold the triangle inequality. Such a distance function called in the literature by semi-metric. (See \cite{tamassy}.) Thus if the set $H$ is
sufficiently smooth, then it can be adopted a metric for it, which arises from the restriction of the Minkowski product to the tangent spaces of $H$.
Let see this more precisely.

The directional derivatives of a function $f:S\longmapsto \mathbb{R}$ with respect to a unit vector $e$ of $S$ can be defined in the usual way, by the existence of the limits for real $\lambda $:
$$
f'_{e}(s)=\lim\limits_{\lambda \mapsto 0}\frac{f(s+\lambda e)-f(s)}{\lambda}.
$$
Let now the generalized Minkowski space is a generalized space-time model, and consider a mapping $f$ on $S$ to $\mathbb{R}$ and the basis ${e_1,\ldots,e_n}$ of Theorem 8. The set of points $F:=\{(s+f(s)e_n)\in V$ \mbox{ for } $s\in S$\} is a so-called {\bf hypersurface} of this space. Tangent vectors of a hypersurface $F$ in a point $p$ are the vectors associated to the directional derivatives of the coordinate functions in the usual way. So
$u$ is a {\bf tangent vector} of the hypersurface $F$ in its point $v=(s+f(s)e_n)$, if it is of the following form
$$
u=\alpha (e+f'_{e}(s)e_n) \mbox{ for real } \alpha \mbox{ and unit vector } e\in S.
$$
The linear hull of the tangent vectors translated into the point $s$ is the tangent space of $F$ in $s$. If the tangent space has dimension $(n-1)$ we call it {\bf tangent hyperplane}.

\begin{lemma}
Let $V$ be a generalized Minkowski space and assume that the s.i.p. $[\cdot,\cdot]|_S$ is continuous. (So the property {\bf s6} holds.) Then the directional derivatives of the real valued function
$$
f:s\longmapsto \sqrt{1+[s,s]},
$$
are
$$
f'_{e}(s)=\frac{\Re{[e,s]}}{\sqrt{1+[s,s]}} \mbox{ for all } i=1,\ldots, n-1.
$$
\end{lemma}
\proof

The considered derivative is:
$$
\frac{f(s+\lambda e)-f(s)}{\lambda}=\frac{\sqrt{1+[s+\lambda e,s+\lambda e]}-\sqrt{1+[s,s]}}{\lambda}=
$$
$$
=\frac{\sqrt{1+[s+\lambda e,s+\lambda e]}\sqrt{1+[s,s]}-(1+[s,s])}{\lambda \sqrt{1+[s,s]}}.
$$
Since $s+\lambda e, s\in S$ and $S$ is a positive subspace thus
$$
0\leq (\sqrt{[s+\lambda e,s+\lambda e]}-\sqrt{[s,s]})^2=
$$
$$
=[s+\lambda e,s+\lambda e]-2\sqrt{[s+\lambda e,s+\lambda e]}\sqrt{[s,s]} +[s,s],
$$
so
$$
[s+\lambda e,s+\lambda e]+[s,s]\geq 2\sqrt{[s+\lambda e,s+\lambda e]}\sqrt{[s,s]}\geq 2|[s+\lambda e,s]|,
$$
and also
$$
[s+\lambda e,s+\lambda e]+[s,s]\geq 2|[s,s+\lambda e]|.
$$
Using these inequalities we get that:
$$
\frac{f(s+\lambda e)-f(s)}{\lambda}\geq \frac{\sqrt{1+2|[s+\lambda e,s]|+|[s+\lambda e,s]|^2}-(1+[s,s])}{\lambda \sqrt{1+[s,s]}}=
$$
$$\frac{1+|[s+\lambda e,s]|-1-[s,s]}{\lambda \sqrt{1+[s,s]}}\geq
\frac{\Re\{[s,s]+\lambda[e,s]\}-[s,s]}{\lambda \sqrt{1+[s,s]}}=\frac{\Re{[e,s]}}{\sqrt{1+[s,s]}}.
$$
But also
$$
\frac{f(s+\lambda e)-f(s)}{\lambda}=
$$
$$
=\frac{(1+[s+\lambda e,s+\lambda e])-\sqrt{1+[s,s]}\sqrt{(1+[s+\lambda e,s+\lambda e])}}{\lambda \sqrt{1+[s+\lambda e,s+\lambda e]}}\leq
$$
$$
\leq \frac{(1+[s+\lambda e,s+\lambda e])-1-|[s,s+\lambda e]|}{\lambda \sqrt{1+[s+\lambda e,s+\lambda e]}}=
$$
$$
=\frac{\Re\{[s+\lambda e,s+\lambda e]\}-|[s,s+\lambda e]|}{\lambda \sqrt{1+[s+\lambda e,s+\lambda e]}}=
$$
$$
=\frac{\Re\{[s,s+\lambda e]+\lambda [e,s+\lambda e]\}-|[s,s+\lambda e]|}{\lambda \sqrt{1+[s+\lambda e,s+\lambda e]}}\leq
$$
$$
\leq\frac{|[s,s+\lambda e]|+\Re\{\lambda [e,s+\lambda e]\}-|[s,s+\lambda e]|}{\lambda \sqrt{1+[s+\lambda e,s+\lambda e]}}=
$$
$$
=\frac{\Re\{[e,s+\lambda e]\}}{\sqrt{1+[s+\lambda e,s+\lambda e]}}.
$$
Now the continuity property {\bf s6} implies that the examined limit exists, and that the differential is
$$
\frac{\Re{[e,s]}}{\sqrt{1+[s,s]}}
$$ as we stated.
\qed

We now apply our investigation to $H$ of a generalized space-time model. As it can be seen easily the explicit form of this hypersurface arises from the above function
$$
f:s\longmapsto \sqrt{1+[s,s]}.
$$
Since its directional derivatives can be determined concretely we can give a connection between the differentiability properties and the orthogonality one.

\begin{lemma}
Let $H$ be the imaginary unit sphere of a generalized space-time model. Then the tangent vectors of the hypersurface $H$ in its point 
$$
v=s+\sqrt{1+[s,s]}e_n
$$ 
form the orthogonal complement $v^{\bot}$ of $v$.
\end{lemma}

\proof
A tangent vector of this space is of the form:
$$
u=\alpha (e+f'_{e}(s)e_n)
$$
where by the previous lemma
$$
f'_{e}(s)=\frac{\Re{[e,s]}}{\sqrt{1+[s,s]}}=\frac{[e,s]}{\sqrt{1+[s,s]}}.
$$
Thus we have:
$$
\left[\alpha \left(e+\frac{[e,s]}{\sqrt{1+[s,s]}}e_n\right),s+t\right]^+=
$$
$$
=\alpha[e,s]+\alpha \left[\frac{[e,s]}{\sqrt{1+[s,s]}}e_n,\sqrt{1+[s,s]}e_n\right]=\alpha([e,s]-[e,s])=0.
$$
So the tangent vectors are orthogonal to the vector $v$. Conversely, if for a vector $u=s'+t'=s'+\lambda e_n$
$$
0=[u,v]=[s',s]+[t',t],
$$
then
$$
[s',s]=-[\lambda e_n,t]=\lambda \sqrt{1+[s,s]}
$$
since $-[t,t]=1+[s,s]$ by the definition of $H$. Introducing the notion
$$
e=\frac{s'}{\sqrt{[s',s']}}
$$
we get that
$$
[e,s]=\left[\frac{s'}{\sqrt{[s',s']}},s\right]=\frac{\lambda}{\sqrt{[s',s']}}\sqrt{1+[s,s]},
$$
implying that
$$
\frac{\lambda}{\sqrt{[s',s']}}=\frac{[e,s]}{\sqrt{1+[s,s]}}=f'_e(s).
$$
In this way
$$
u=\sqrt{[s',s']}\left(\frac{s'}{\sqrt{[s',s']}}+\frac{\lambda}{\sqrt{[s',s']}}e_n\right)=\alpha(e+f'_e(s)e_n).
$$
This last equality shows that a vector of the orthogonal complement is a tangent vector as we stated.
\qed

We define now the Finsler space type structure for a hypersurface of a generalized space-time model.

\begin{defi}

If $F$ is a hypersurface of a generalized space-time model for which the following properties hold:
\begin{description}
\item i, in every point $v$ of $F$, there is an (unique) tangent hyperplane $T_v$ for which the restriction of the Minkowski product
$[\cdot,\cdot]^+_v$ is positive,

\item ii, the function $ds^2_v:=[\cdot,\cdot]^+_v:F\times T_v\times T_v\longrightarrow \mathbb{R^+}$
$$
ds^2_v:(v,u_1,u_2)\longmapsto [u_1,u_2]^+_v
$$
varying differentiable with the vectors $v\in F$ and $u_1,u_2\in T_v$,
\end{description}
then we say that the pair $(F, ds^2)$ is a Minkowski-Finsler space with semi-metric $ds^2$ embedding into the generalized space-time model $V$.
\end{defi}

Naturally "varying differentiable with the vectors $v,u_1,u_2$" means that for every $v\in T$ and pairs of vectors $u_1,u_2\in T_v$ the function $[u_1,u_2]_v$ is a differentiable function on $F$.

\begin{theorem}
Let $V$ be a generalized space-time model. Let $S$ be a continuously differentiable s.i.p. space then $(H^+,ds^2)$ is a Minkowski-Finsler space.
\end{theorem}

\proof If the s.i.p. of $S$ is a continuously differentiable one, then the norm is differentiable twice (See Theorem 2.). This also implies the
continuity of the s.i.p. and so by Lemma 4 we know that there is an unique tangent hyperplane at each point of $H$. By Theorem 10 we get that the
Minkowski product restricted to a tangent hyperplane is positive so the first assumption of the definition is valid.

To prove the second condition consider the product: $ [u_1,u_2]^+_v, $ where $v$ is a point of $H$ and $u_1$,$u_2$ are two vectors on its tangent
hyperplane. Then by Lemma 4 we have:
$$
u_i=\alpha_i\left(s_i+\frac{[s_i,s_v]}{\sqrt{1+[s_v,s_v]}}e_n\right), \mbox{ for } i=1,2.
$$
Here the vectors $s_1,s_2,s_v$ are in $S$ and $v=s_v+\sqrt{1+[s_v,s_v]}e_n$. Thus the examined product is
$$
[u_1,u_2]^+_v=\alpha_1\alpha_2\frac{[s_1,s_2](1+[s_v,s_v])-[s_1,s_v][s_2,s_v]}{(1+[s_v,s_v])}.
$$
Since the function
$$
[s_v,s_v]=([v,e_n]^+)^2-1
$$
is a continuously differentiable function of $v$, and $[s_1,s_2]$ ( by our assumption) is also continuously differentiable of its arguments, we have to prove only, that the map sending $u_i$ to $s_i$ also holds this property. But this latter fact is a consequence of the observation that the map $u\mapsto s$ is a projection so it is linear.
\qed

\subsection{ The geometry of $H^+$. }

Our next goal will be to give a characterization of the isometries of the Minkowski-Finsler manifold $H^+$. For these we need some further
definitions. The following concept of linear isometry in any generalized Minkowski space is usable.

\begin{defi}
A linear isometry $f:H^+\longrightarrow H^+$ of $H^+$ is the restriction to $H^+$ of a linear map $F:V\longrightarrow V$ which preserves the
Minkowski product and which sends $H^+$ onto itself.
\end{defi}

We note that in this definition a linear mapping $F$  restricted to $S$ gives an isometry between $S$ and its image $F(S)$ implying that this image is a normed space with respect to those s.i.p. which raised from the s.i.p. of $S$. This isometry is stronger than the usual one, in which we need only the equality of the norm of the corresponding vectors.  As we can see in the paper of Koehler:
\begin{theorem}[\cite{koehler}]
In a smooth Banach space a mapping is an isometry if and only if it preserves the (unique) s.i.p..
\end{theorem}
Thus if the norm is at least smooth then the two sense of linear isometry are coincide. Also Koehler proved that if the generalized Riesz-Fischer
representation theorem is valid in a normed space then for every bounded linear operator $A$ has a generalized adjoint $A^T$ defined by the equality:
$$
[A(x),y]=[x,A^T(y)] \mbox{ for all } x,y \in V.
$$
This mapping is the usual Hilbert space adjoint if the space is an i.p. one. In this more general setting this map is not usually linear but it still has some interesting properties. The assumption for the s.i.p. in Koehler paper \cite{koehler} is that the space should be a smooth and uniformly convex Banach space. It is well known that uniform convexity implies strict convexity. On the other hand, we now also take into consideration (see \cite{wilansky} p. 111) that every strictly convex finite-dimensional normed vector space is uniformly convex so for the rest of the section we shall assume that the normed space $S$ with respect to its s.i.p. is strictly convex and smooth. It is convenient to characterize strict convexity of the norm in terms of s.i.p. properties. E.Berkson \cite{berkson} states, and it can be proved simply, that:

\begin{lemma}[\cite{berkson}]
An s.i.p. space is strictly convex if and only if whenever $[x,y]=\|x\| \|y\|$ where $x,y\neq 0$, then $y=\lambda x$ for some real $\lambda>0$.
\end{lemma}

Now we prove the following theorem:

\begin{theorem}
Let V be a generalized space-time model. Assume that the subspace $S$ is a strictly convex, smooth normed space with respect to the norm arisen from
the s.i.i.p.. Then the s.i.p. space $\{V,[\cdot,\cdot]^-\}$ is also smooth and strictly convex. Let $F^T$ be the generalized adjoint of the linear
mapping $F$ with respect to the s.i.p. space $\{V,[\cdot,\cdot]^-\}$, and define the idempotent linear mapping $J:V\longrightarrow V$  by the
equalities $J|S=id|S$, $J|T=-id|_T$. The map $F|_H=f:H\longrightarrow H$ is a linear isometry of the upper sheet $H^+$ of $H$ if and only if it is
invertible, satisfies the equality:
$$
F^{-1}=JF^{T}J,
$$
moreover takes $e_n$ into a point of $H^+$.
\end{theorem}

\proof First we prove that the embedding normed space $\{V,
[\cdot,\cdot]^-\}$ is also smooth and strictly convex. The equality
$1=[s+t,s+t]^-=[s,s]-[t,t]=[s,s]+\|t\|^2$ shows that the unit balls
of the two norms are smooth at the same time. To prove strict
convexity consider
$$
[s+t,s'+t']^-=\|s+t\|^-\|s'+t'\|^-.
$$
Since $\dim T=1$, we can assume that $t'=\lambda t$ for some real $\lambda $. Thus we get the equality:
$$
[s,s][s',s']=[s,s']^2+[t,t]([s',s']-2\lambda [s,s'] +\lambda ^2[s,s]).
$$
By Cauchy-Schwartz inequality we have:
$$
[s',s']-2\lambda [s,s'] +\lambda ^2[s,s]\geq \left(\sqrt{[\lambda s,\lambda s]}-\sqrt{[s',s']}\right)^2\geq 0,
$$
so
$$
0\leq [s,s']^2\leq [s,s][s',s']=[s,s']^2+[t,t]([s',s']-2\lambda [s,s'] +\lambda ^2[s,s])\leq [s,s']^2
$$
implying that
$$
[t,t]([s',s']-2\lambda [s,s'] +\lambda ^2[s,s])=0.
$$
If $[t,t]=0$ then $t=t'=0$ and from the strict convexity of $S$ we get that there is a real $\mu >0$ with $s'=\mu s$. For this $\mu $ we have $s'+t'=\mu (s+t)$, too. So we can assume that $[t,t]\neq 0$ and thus
$$
[s,s][s',s']=[s,s']^2 \mbox{ and }  [s',s']-2\lambda [s,s'] +\lambda ^2[s,s])=0
$$
hold paralelly. But $S$ is a strictly convex space so for a nonzero $s$ there is a real $\mu >0$ with $s'=\mu s$, again. But this also implies
$$
0=(\mu -\lambda)^2[s,s],
$$
showing that $\mu=\lambda $ and $s'+t'=\mu (s+t)$. Using Lemma 5, we get the strict convexity
of the embedding normed space.

Let $F$ be a linear isometry of $H$. It is clear that the linear operator $J$ transforms the Minkowski product into the s.i.p. of the embedding
space. Precisely we have:
$$
[v,w]^+=[v,Jw]^-.
$$
Now using the existence of the adjoint operator, the following calculation:
$$
[v,Jw]^-=[v,w]^+=[Fv,Fw]^+=[Fv,JFw]^-=[v,F^TJFw]^-
$$
holds for each pair of vectors $v$ and $w$. But the embedding space is a nondegenerate one, thus we get the equality:
$$
J=F^TJF \mbox{ or equivalently }
$$
$$
F^{-1}=JF^{T}J.
$$
By its definition the last condition on $F$ also holds.

Conversely, if $F$ is a linear mapping satisfying the condition of the theorem then preserves the Minkowski product. In fact,
$$
[Fv,Fw]^+=[Fv,JFw]^-=[v,F^TJFw]^-=[v,Jw]^-=[v,w]^+.
$$
It takes the hyperboloid $H$ homeomorphically onto itself implying that it takes a sheet onto a sheet. Our last condition guarantees that $F(H^+)=H^+$ and $F$ is a linear isometry of $H^+$ as we stated.
\qed

As it can be seen from the formula of Theorem 13 the generalized adjoint of a linear isometry is a linear transformation. We also note that Theorem
13 in the i.p. case gives the characterization of the isometries of the hyperbolic space of dimension $(n-1)$.

It is not clear that there is or is not a non-pseudo Euclidean generalized Minkowski space for which the group of linear isometries acts transitively
on $H^+$. But if the answer is yes and so the Minkowski-Finsler geometry of $H^+$ is linearly homogeneous, then we can compute the Minkowski-Finsler
distance. Now we determine the distance function $d:H^+\times H^+\longrightarrow \mathbb{R^+}$ of a linearly homogeneous Minkowski-Finsler space
$H^+$.

Before the calculation we recall some known concept on classical Finsler spaces. We assume that the s.i.i.p. restricted into $S$ is continuously
differentiable. In a connected Finsler space any point has a distance from any other point of the space (see e.g. \cite{tamassy}). By our terminology
the distance can be got in the following analogous way.

\begin{defi} Denote by $p,q$ a pair of points in $H^+$ and consider the set $\Gamma _{p,q}$ of equally oriented piecewise differentiable curves $c(t)$
$a\leq t\leq b$ of $H^+$ emanating from $p$ and terminating at $q$. Then the Minkowskian-Finsler distance of these points is:
$$
\rho (p,q)=\inf \left\{ \int \limits _{a}^{b}\sqrt{[\dot{c}(x),\dot{c}(x)]^+_{c(x)}}dx \mbox{ for } c\in \Gamma _{p,q}\right \},
$$
where $\dot{c}(x)$ means the tangent vector of the curve $c$ in its point $c(x)$.
\end{defi}

We would like to examine the influence of a linear isometry to the Minkowski-Finsler distance. It is easy to see that this distance holds the
triangle inequality thus it is a metric on $H^+$. (See \cite{tamassy}.)

\begin{defi}
A topological isometry $f:H\longrightarrow H$ of $H$ is a homeomorphism of $H$ which preserves the Minkowski-Finsler distance between each pair of
points of $H$.
\end{defi}

First we reformulate the length of a path as follows. The Minkowski-Finsler semi-metric on $H^+$ is the function $ds^2$ which assigns at each point
$v\in H^+$ the Minkowski product which is the restriction of the Minkowski product to the tangent space $T_v$. This positive Minkowski product varies
differentiable with $v$. Let $U\leq V$ be a subspace and consider a map $f:U\longrightarrow V$. If it is a totally differentiable map (with respect
to the norm of the embedding $n$-space in the sense of Frechet) then $f(T_v)=T_{f(v)}$ for the tangent spaces at $v$ and $f(v)$, respectively and one
can define the pullback semi-metric $f^\star(ds^2)$ at the point $v$ by the following formula:
$$
f^\star(ds^2)_v(u_1,u_2)=ds^2_{f(v)}(Df(u_1),Df(u_2))=[Df(u_1),Df(u_2)]^+_{f(v)}.
$$
The square root $ds$ of the semi-metric function defined by $\sqrt{ds^2_v(u,u)}$ is the so called length element and the length of a path is the
integral of the pullback length element by the differentiable map $c:\mathbb{R}\longrightarrow V$. This implies that if a linear isometry leaves
invariant the Minkowski-Finsler semi-metric by the pullback then it preserves the integrand and thus preserves the integral, as well. Let now $F$ be
a linear isomorphism and its restriction to $H^+$ is $f$. Compute the pullback metric as follows:
$$
f^\star(ds^2)_v(u_1,u_2)=ds^2_{f(v)}(Df(u_1),Df(u_2))=[Df(u_1),Df(u_2)]^+_{f(v)}=
$$
$$
=[DF(u_1),DF(u_2)]^+_{F(v)}=[F(u_1),F(u_2)]^+_{F(v)}
$$
because $F$ is linear. But it preserves the Minkowski product therefore we conclude that
$$
[F(u_1),F(u_2)]^+_{F(v)}=[u_1,u_2]^+_{v}=(ds^2)_v(u_1,u_2).
$$
This proves the following theorem:
\begin{theorem}
A linear isometry of $H^+$ is a topological isometry on it, too.
\end{theorem}
In the proof of this theorem we also proved that a linear isometry is a Finsler isometry, in the sense that it is a diffeomorphism of $H$ onto $H$
which preserves the Minkowski-Finsler metric function. In a Riemann space the two kind of isometries (the topological and Riemannian one) are
equivalent. This is the result of Myers and Steenrod (See in \cite{myers}). The analogous theorem on Finsler spaces was proved by Deng and Hou in
\cite{deng}. This latter one states that the two concepts of isometry are equivalent for a Finsler space, too.

In the following theorem we impose the condition of linear homogeneity of $H^+$. Thus we state:
\begin{theorem}
Let $V$  be a generalized space-time model. Consider that the normed space $S$ is strictly convex and smooth and the group of linear isometries of
$H^+$ acts transitively on $H^+$. Let denote the Minkowski-Finsler distance of $H^+$ by $d(\cdot,\cdot)$. Then the following statement is true:
$$
[a,b]^+=-ch(d(a,b)) \mbox{ for } a,b\in H^+.
$$
\end{theorem}
\proof In a Finsler space a function preserving the distance function transforms geodesics to geodesics. (See in \cite{bao}.) In our case this is
also true since this fact basically determined by the definition of the distance and the smoothness properties which are same in both cases.
Since our space is homogeneous and linear isometry preserves the distance by Theorem 14, we can assume that $a=e_n$. Let now $b\neq a$ and consider
the 2-plane $<a,b>$ spanned by the vectors $a$ and $b$. The restriction of the s.i.i.p. to the plane $<a,b>$ is an i.i.p. thus the restricted Finsler
function is a Riemannian one. So the intersection $H\cap<a,b>$ is a hyperbole in the embedding Euclidean two space, thus we can parameterize the
points of a path from $a$ to $b$ by
$$
c(t)=sh(\tau)e+ch(t)e_n \mbox{ for } t\in [0,1],
$$
with $c(0)=a$ and $c(1)=b$. The length of an arc from $0$ to $x$ is:
$$
\int \limits _0^x\sqrt{ch^2(\tau)-sh^2(\tau)}d\tau =x,
$$
showing that the points of this arc satisfy the triangle inequality by equality. Consequently it is a geodesic on $H^+$ therefore its arc-length is the distance of the point $a$ and $c(x)$. On the other hand we also have:
$$
[a,b]^+=[e_n,sh(1)e+ch(1)e_n]^+=[e_n,ch(1)e_n]=-ch(1)=
$$
$$
=-ch(d(a,c(1))=-ch(d(a,b)).
$$
\qed

\begin{center}
\'Akos G.Horv\'ath,\\
 Department of Geometry \\
Budapest University of Technology and Economics\\
1521 Budapest, Hungary
\\e-mail: ghorvath@math.bme.hu
\end{center}


\begin{thebibliography}{99}

\bibitem{alonso1} Alonso, J., Benitez, C.:  {\em Orthogonality in normed linear spaces: a survey. Part I. Main properties.} Extracta Math. {\bf 3} (1988), 1--15.

\bibitem{alonso2} Alonso, J., Benitez, C.:  {\em Orthogonality in normed linear spaces: a survey. Part II. relation between main orthogonalities.} Extracta Math. {\bf 4} (1989), 121--131.

\bibitem{bao} Bao, D., Chern S.S., Shen Z.: {\em An introduction to Riemannian-Finsler Geometry} Springer-Verlag, Berlin, 1999.

\bibitem {berkson} Berkson, E.: {\em Some type of Banach spaces, Hermitian operators and Bade functionals.} Trans. Amer. Math. Soc. {\bf 116} (1965), 376--385.

\bibitem{cannon} Cannon, J.W., Floyd W.J., Kenyon, R., Parry, W.R.: {\em Hyperbolic geometry.} http://citeseerx.ist.psu.edu/ viewdoc/ summary? doi=10.1.1.31.1601.

\bibitem{clarkson} J.A.Clarkson, {\em Uniformly convex spaces.} Trans. Amer. Math. Soc. {\bf 40} (1936), 396--414.

\bibitem{day} Day, M.M.: {\em Polygons circumscribed about closed convex curves} Trans. Amer. Math. Soc. {\bf 62} (1947), 315--319.

\bibitem {deng} Deng, S., Hou, Z.: {\em The group of isometries of a Finsler space} Pacific J. of Math. {\bf 207/1} (2002), 149--155.

\bibitem {diminnie} Diminnie, C.R.: {\em A New Orthogonality Relation for Normed Linear Spaces.} Math. Nachr.{\bf 114} (1983) 197--203.

\bibitem {gahler} Gahler, S.: {\em Lineare 2-normierte Raume.} Math. Nachr. {\bf
28} (1964) 1--43.

\bibitem{giles} Giles, J.R.: {\em Classes of semi-inner-product spaces.} Trans. Amer. Math.Soc. {\bf 129/3} (1967), 436--446.

\bibitem{gohberg} Gohberg, I.,Lancester, P., Rodman, L.: {\em Indefinite Linear Algebra and Applications} Birkhauser, Basel-Boston-Berlin 2005.

\bibitem{james} James, R.C.: {\em Orthogonality in normed linear space.}
Duke Math. J. {\bf 12} (1945), 291--301.

\bibitem {knowles} Knowles R.J., Cook, T.A.: {\em Some results on Auerbach bases for finite dimensional normed spaces.} Bull. Soc. Roy. Sci. Li\`{e}ge {\bf 42} (1973), 518--522.

\bibitem {koehler} Koehler D.O.: {\em A Note on Some Operator Theory in Certain Semi-Inner-Product Spaces.} Proc. Amer. Math. Soc. {\bf 30(2)} (1971) 363--366.

\bibitem {l-g} Gruber P.M.- Lekkerkerker C.C.: {\em Geometry of num\-bers.\/ } \hfill \break North-Holland Amsterdam-New York-Oxford-Tokyo 1987.

\bibitem {lenz} Lenz, H.: {\em Eine Kennzeichnung des Ellipsoids.} Arch. Math. {\bf 8} (1957), 209--211.

\bibitem {lumer} Lumer, G.:{\em Semi-inner product spaces} Trans. Amer. Math. Soc. {\bf 100} (1961), 29-43.

\bibitem {lumer2} Lumer, G.:{\em On the isometries of reflexive Orlicz spaces.} Ann. Inst. Fourier, Grenoble {\bf 13} (1963) 99--109.

\bibitem{martini}  Martini, H.: {\em Shadow boundaries of convex bodies.} Discrete
Math. {\bf 155}, (1996) 161-172.

\bibitem {martini-swanepoel 1} Martini, H., Swanepoel, K., Weiss, G.: {\em The
geometry of Minkowski spaces - a survey. Part I.} Expositiones
Mathematicae {\bf 19} (2001), 97-142.

\bibitem {martini-swanepoel 2} Martini, H., Swanepoel, K.: {\em The Geometry of Minkowski
Spaces - A survey. Part II.} Expositiones Mathematicae {\bf 22(2)}
(2004), 93-144.

\bibitem {mcshane} McShane, E.J. {\em Linear functionals on certain Banach spaces.} Proc. Amer.Math.Soc., Vol. {\bf 1} (1950), 402--408.

\bibitem {milicic} Milicic, P.M.: {\em Sur le q-angle daus um espace norme.} Mat. Vesnik {\bf
45} (1993) 43--48.

\bibitem {myers} Myers, S.B., Steenrod, N.: {\em The group of isometries of a Riemannian manifold} Ann. of Math., {\bf 40} (1939), 400-416.

\bibitem {nath} Nath, B.: {\em On a generalization of semi-inner product spaces} Math. J. Okayama Univ. {\bf 15/1} (1971), 1--6.

\bibitem {partington} Partington, J.R.: {\em Orthogonality in normed spaces.} Bull. Austral Math. Soc. {\bf 33} (1986), 449--455.

\bibitem {shoja} Shoja , Maeheri, H.: {\em General Orthogonality in Banach spaces.} Int. J. Math. Analysis, {\bf 1(12)} (2007) 553-556.

\bibitem {tamassy} Tam\'assy, L.: {\em Finsler spaces corresponding to distance spaces} Proc. of the Conf., Contemporary geometry and related Topics, Belgrade, Serbia and Montenegro, June 26--July 2, (2005), 485--495.

\bibitem {taylor} Taylor, A.E.: {\em A geometric theorem and its application to biorthogonal systems.} Bull. Amer. Math. Soc. {\bf 53} (1947), 614--616.

\bibitem {wilansky} Wilansky, A.: Functional analysis. {\em Blaisdell, New Yourk} 1964.





\end{thebibliography}
\end{document}